\documentclass[a4paper]{article}
\usepackage[utf8]{inputenc}
\usepackage[T1]{fontenc}
\usepackage[english]{babel}
\usepackage[pdfusetitle,hidelinks]{hyperref}
\usepackage{amsmath,amsfonts,amssymb,amsthm}
\usepackage{graphicx}
\usepackage[noabbrev,capitalize]{cleveref}
\usepackage{tikz-cd}
\usepackage[shortlabels]{enumitem}
\usepackage{macros}
\usepackage{authblk}
\usepackage{theorems}

\crefname{enumi}{}{}
\setlist[itemize]{itemsep=0pt, label={---}}
\newcommand{\citep}{\cite}

\usepackage{titlesec}
\titleformat{\paragraph}[runin]{\normalfont\bf}{\thesubsection}{1.5ex}{}[.]
\crefformat{equation}{(#2#1#3)}

\title{Hypercubical manifolds in homotopy type theory}
\author[1]{Samuel Mimram}
\author[1]{Émile Oleon}
\affil[1]{LIX, CNRS, École polytechnique, Institut Polytechnique de Paris, Palaiseau, France}
\date{}

\begin{document}
\maketitle

\begin{abstract}
  Homotopy type theory provides a logical framework in which geometric constructions and proofs can be carried out synthetically: in this setting, types correspond to spaces up to homotopy, and proofs to homotopy-invariant constructions. Within this context, we introduce a type corresponding to the \emph{hypercubical manifold}, a space first described by Poincaré in 1895. This manifold is interesting because it offers an approximation of the quaternion group $Q$, in the sense that it represents the first step toward the construction of a cellular resolution of $Q$. To validate our definition, we show that it satisfies the expected property: it is the homotopy quotient of the 3-sphere $\S3$ under the natural action of $Q$. Establishing this result is non-trivial, requiring subtle combinatorial computations based on the flattening lemma, thereby illustrating the constructive power of homotopy type theory. Finally, extending this construction, we introduce higher-dimensional generalizations of the manifold, which provide increasingly precise cellular approximations of $Q$, and converge toward a delooping of $Q$.

\end{abstract}


\section{Introduction}
Homology is a family of abelian groups, defined for spaces such as topological spaces or simplicial sets, that remain invariant under homotopy equivalence, \ie under continuous deformation~\citep{maclane2012homology}. After its introduction, the question arose of how precise this invariant is and, in particular, whether it could uniquely characterize spheres~\citep{poincare1900second}. The answer is negative: there are spaces, known as \emph{homology spheres}, that share the same homology groups as spheres, but are not homotopy equivalent to them.
Constructing such examples led Poincaré to investigate a family of spaces obtained by gluing faces of Platonic solids~\citep{poincare1895analysis}, which includes the homology sphere~\citep{poincare1904cinquieme}.

\paragraph{The hypercubical manifold}
One such space, the \emph{hypercubical manifold}~$\HM$, can be described as

\begin{enumerate}[(A)]
\item\label{hm-def-cw} a 3-dimensional cube in which each pair of opposite square faces is identified after a quarter-turn rotation.
\end{enumerate}

\noindent
Modern expositions of this space can be found in~\citep{hypercubique,sarem}. Its fundamental group is the quaternion group~$Q$ and its universal cover is the 3-sphere $\S3$. It is well known that any (sufficiently nice) space can be recovered as the quotient of its universal cover by the canonical action of its fundamental group~\cite[Proposition 1.40]{hatcher}. Consequently, we have an alternative description of~$\HM$ as

\begin{enumerate}[(A),resume]
\item\label{hm-def-quotient} the quotient of the 3-sphere $\S3$ by a free action of the quaternion group~$Q$.
\end{enumerate}

\paragraph{Defining the hypercubical manifold in homotopy type theory}
In recent years, a variant of dependent type theory, \emph{homotopy type theory}, has been introduced~\citep{hottbook}. In this setting, types can be interpreted as homotopy types, \ie spaces up to homotopy equivalence, thus allowing one to perform constructions on spaces in a synthetic way. Alternative versions of this type theory, such as the one implemented in Cubical Agda~\citep{vezzosi2021cubical}, allow for a generic implementation of higher inductive types (HITs), where one can use higher constructors in the definition of types, corresponding to attaching cells as in CW-complexes. Our goal in this paper is to define the hypercubical manifold in homotopy type theory (as well as higher-dimensional variants) and show that it satisfies the expected properties: in particular, we will show the equivalence of the above two definitions within the theory.

The hypercubical manifold can easily be constructed as a HIT~\citep{dylan} by translating the standard description \cref{hm-def-cw} of the space as a CW-complex (many types corresponding to traditional spaces can be defined in this way~\cite[Section 2.4]{vezzosi2021cubical}).
%
However, constructing $\HM$ as a quotient of the 3-sphere by an action of $Q$, following the second definition \cref{hm-def-quotient}, seems difficult to achieve, because there is no direct way of defining such an action. In order to do so, we would need to have a \emph{coherent} description of it, that is as a map~$\phi:\B Q\to\U$ where $\B Q$ is a type internally representing the group $Q$ (this is called a \emph{delooping} of~$Q$) and~$\U$ is the universe of all types, whose elements can be thought of here as spaces. Classical definitions of $\B Q$ (as a type of $Q$-torsors~\citep{warn2023eilenberg}, or as an Eilenberg-MacLane space~\citep{licata2014eilenberg}) have a universal property that only allows one to eliminate into groupoids, and there is thus no direct way of defining the map~$\phi$ because $\S3$ is not $n$-truncated, for any~$n$.


Here, we explain how we can still manage to define the coherent action of $Q$ on $\S3$ in homotopy type theory, from a map $\HM\to\B Q$ whose fiber is $\S3$, and we show that the resulting homotopy quotient is equivalent to the definition of $\HM$ as a HIT. This implies that the latter definition is ``correct'': it is a quotient of $\S3$ by $Q$ as expected. This construction uses two main ingredients: the action-fibration duality (\cref{action-fibration}) in order to construct the map~$\phi$, and the flattening lemma (\cref{flattening-pushout}) in order to compute its fibers. Note that Definition~\ref{hm-def-quotient} uses a strict quotient, whereas we only have access to homotopy quotients in homotopy type theory, but the two coincide in the case of free actions, which can be seen through the Borel construction~\cite[Section 2.2]{adem2002topics}.

Previous results imply that the groupoid truncation of $\HM$ is a model of $\B Q$ (this is a variant of the traditional construction of Eilenberg-MacLane spaces in homotopy type theory by~\cite{licata2014eilenberg}). However, the truncation is formal and thus, for the same reasons as above, it is unsatisfactory from a computational point of view. From this perspective, it is much better to use the sequence of spaces obtained from the join construction in homotopy type theory~\cite{rijke2017join,rijke2018classifying}, which provide increasingly good approximations of $\B Q$ and actually converge toward~$\B Q$.
If we apply the join construction starting from~$\HM$, we thus obtain a sequence of types~$\HM[n]$, which provide increasingly good approximations of~$\B Q$ in the sense that the canonical map $\HM[n]\to\B Q$ is $(4n-2)$-connected. We observe here that $\HM[n]$ comes from a quotient of the sphere $\S{4n-1}$ by a free action of~$Q$ and, from this point of view, can thus be thought of as a higher-dimensional analogue of the hypercubical manifold.
The types $\HM[n]$ are interesting because they provide a good description of $\B Q$ as a HIT which is non-recursive (in particular, it avoids resorting to formal groupoid truncation), relatively ``small'' in terms of the number of constructors in each dimension, and expected to be cellular and thus more amenable to computations (such as cohomology groups, as explained in the conclusion of the article).

%

\paragraph{Plan of the paper}
We begin by recalling the traditional definition of the hypercubical manifold in topological spaces in \cref{topological-definition}, as well as basic notation and constructions in homotopy type theory in \cref{hott}. The main contributions of the paper can be found in \cref{hypercubical-hott}, where we define the hypercubical manifold in homotopy type theory, and show that the fiber of the canonical map $\HM\to\B Q$ is $\S3$ by iteratively using the flattening lemma on skeleta of $\HM$. We then deduce all expected properties from this fiber sequence in \cref{higher}. Finally, we introduce and motivate higher-dimensional generalizations of the hypercubical manifold and conclude in \cref{sec:conclusion}.

\section{Topological definition of the hypercubical manifold}
\label{topological-definition}
Before turning to our implementation of the hypercubical manifold in homotopy type theory, we first recall its classical construction, based on a free action of the quaternion group on the 3-sphere. This whole section is to be read in the setting of classical set-theoretic mathematics (as opposed to homotopy type theory): all quotients are taken to be strict ones, etc.

\subsection{The quaternion group}
Consider the real 4-dimensional vector space~$\H$ generated by the basis vectors $1$, $i$, $j$, and $k$. This space is an algebra under the Hamilton product, with $1$ as neutral element, and given on basis elements by $i^2=j^2=k^2=-1$, $ij=k$, $ji=-k$, $jk=i$, $kj=-i$, $ki=j$ and $ik=-j$. The \emph{quaternion group} is the 8-element group consisting of $1,i,j,k$ and their negatives. This group can easily be shown to admit the presentation
\begin{equation}
  \label{q-std-pres}
  Q=\pres{e,i,j,k}{i^2=e,j^2=e,k^2=e,ijk=e,e^2=1}
\end{equation}
where $e$ denotes~$-1$. The following alternative presentation will prove useful later in the paper.

\begin{lemma}
  \label[lemma]{quaternion-pres}
  The quaternion group admits the presentation
  \begin{equation}
    \label{q-pres}
    Q=\pres{i,j}{i=jij,j=iji}
  \end{equation}
\end{lemma}
\begin{proof}
  We first show that the group admits the presentation
  \begin{equation}
    \label{q-pres2}
    Q=\pres{i,j}{i^4=1,i^2=j^2,iji=j}
  \end{equation}
  We have the following equivalent presentations
  \begin{align*}
    Q&=\pres{e,i,j,k}{i^2=e,j^2=e,k^2=e,ijk=e,e^2=1}
    \\
    Q&=\pres{i,j,k}{j^2=i^2,k^2=i^2,ijk=i^2,i^4=1}&&\text{using $e=i^2$}
    \\
    \intertext{
      In the second presentation above, $ij=k$ is derivable because we have $ijk=i^2$, multiplying by~$k^3$ yields $ijk^4=i^2k^3$, and therefore $ij=k$ since $k^4=i^4=1$ and $i^2k^2=i^4=1$. We thus have the following presentations:
      }
    Q&=\pres{i,j,k}{j^2=i^2,k^2=i^2,ijk=i^2,i^4=1,ij=k}
    \\
    Q&=\pres{i,j}{j^2=i^2,(ij)^2=i^2,ijij=i^2,i^4=1}&&\text{using $k=ij$}
    \\
    \intertext{
      The relation $(ij)^2=i^2$ is now redundant, so that we can remove it, and then, from the relation $ijij=i^2=j^2$, we can derive the equivalent relation $iji=j$ by dividing by $j$, and we obtain \cref{q-pres2}:
    }
    Q&=\pres{i,j}{i^4=1,i^2=j^2,iji=j}&&\text{\qquad\qquad\qquad\qquad\qquad\qquad\qquad}
  \end{align*}

  We can now show that the presentation \cref{q-pres2} is equivalent to \cref{q-pres}. From \cref{q-pres2}, we can derive the missing relation $i=jij$
  by
  \begin{align*}
    jij&=iji^2j&&\text{by $j=iji$}\\
    &=ij^4&&\text{by $i^2=j^2$}\\
    &=i&&\text{by $j^4=i^4=1$}
  \end{align*}
  Conversely, from \cref{q-pres} we can derive the relations of \cref{q-pres2} since $i^2=jiji=j^2$, and $j=iji=jij^2i=ji^4$ which implies $1=i^4$ by division by~$j$.
\end{proof}

\subsection{The hypercubical manifold}
\label{hypercubical-manifold}
We write $C=\setof{t1+xi+yj+zk}{t,x,y,z\in[-1,1]}\subseteq\H$ for the \emph{unit cube} in $\H$, and $\partial C\subseteq C$ for its boundary (consisting of the points such that at least one of the coordinates is either $-1$ or $1$). Note that $\partial C$ is homotopy equivalent to the 3-sphere $\S3$. There is a natural action of~$Q$ on $C$ by left multiplication. This action is not free since it fixes $0$, but it becomes free when restricted to~$\partial C$.
We write
\[
  \HM\qdefd\partial C/Q
\]
for the quotient space of~$\partial C\simeq\S3$ under this action, and refer to it as the \emph{hypercubical manifold}.
It is not immediately clear from this description that $\HM$ is actually a manifold (as opposed to a topological space), but this will not play a role here.
Because the action is free, and properly discontinuous as~$Q$ is finite, the quotient map $\partial C\to\HM$ is covering with fiber $Q$~\cite[Proposition 1.40]{hatcher}, \ie we have a fiber sequence
\[
  \begin{tikzcd}
    Q\ar[r]&\partial C\ar[r]&\HM
  \end{tikzcd}
\]
from which we deduce that the fundamental group of $\HM$ is $Q$~\cite[Theorem 4.41]{hatcher}.

\label{hm-fundamental-group}

The original definition of~$\HM$~\cite[p. 66, section 13, troisième exemple]{poincare1895analysis} is quite different in nature, with the equivalence explained in~\citep{hypercubique}. Indeed, the hypercubical manifold can alternatively be described as a 3-dimensional cube, whose opposite faces are identified after a quarter-turn rotation:
\begin{equation}
  \label{hm}
  \begin{tikzcd}
    &\ar[dl,"w"']a\ar[dd,dotted,"x",near start]\ar[rr,"y"]&&b\\
    b&&\ar[ll,"z",near start]a\ar[dd,"w",near start]\ar[ur,"x"]\\
    &b&&\ar[ll,dotted,"w",near end]\ar[dl,"y"]a\ar[uu,"z"']\\
    a\ar[uu,"y"]\ar[ur,dotted,"z"]\ar[rr,"x"']&&b
  \end{tikzcd}
\end{equation}
Hence $\HM$ can be represented as a CW-complex with two 0-cells ($a$ and $b$), four 1-cells ($x$, $y$, $z$ and~$w$), three 2-cells ($\alpha$, $\beta$ and $\gamma$, corresponding to the front, right and top faces)
\begin{equation}
  \label{2cells}
  \begin{tikzcd}
    b\ar[dr,phantom,"\alpha"]&\ar[l,"z"']a\ar[d,"w"]\\
    a\ar[u,"y"]\ar[r,"x"']&b
  \end{tikzcd}
  \qquad\qquad\qquad
  \begin{tikzcd}
    b\ar[dr,phantom,"\beta"]&\ar[l,"w"']a\ar[d,"x"]\\
    a\ar[u,"y"]\ar[r,"z"']&b
  \end{tikzcd}
  \qquad\qquad\qquad
  \begin{tikzcd}
    b\ar[dr,phantom,"\gamma"]&\ar[l,"w"']a\ar[d,"y"]\\
    a\ar[u,"z"]\ar[r,"x"']&b
  \end{tikzcd}
\end{equation}
and one 3-cell ($A$, representing the interior of the cube).
This definition allows for a direct computation of the fundamental group $\pi_1(\HM)$. We first contract $w$ to an identity (thus identifying $a$ and~$b$) to obtain the presentation
\[
  \pi_1(\HM)
  \qeq
  \pres{x,y,z}{x=zy,y=xz,z=yx}
\]
where the three relations correspond to the 2-cells \cref{2cells}. By removing $z$ (which is equal to $yx$), renaming $x$ to~$i$ and $y$ to~$j$, we recover the presentation of \cref{quaternion-pres} of the quaternion group.

\subsection{Higher hypercubical manifolds}
\label{higher-topological-hypercubical}
We now briefly discuss higher-dimensional generalizations of the hypercubical manifold. These arise quite naturally from our work, but to the best of our knowledge, have not been previously considered in the literature. Writing $\HM[1]$ for the hypercubical manifold, we observe that $\HM[1]$ is a ``good approximation'' of~$Q$ up to dimension~$2$, in the sense that we have
\begin{align*}
  \pi_0(\HM[1])=\pi_0(\B Q)&=1
  &
  \pi_1(\HM[1])=\pi_1(\B Q)&=Q
  &
  \pi_2(\HM[1])=\pi_2(\B Q)&=1
\end{align*}
More generally, we would like to have cellular spaces $\HM[n]$ which coincide with~$Q$ up to dimension~$f(n)$ where $f$ is a strictly increasing function on the natural numbers, thus constructing a resolution of~$Q$. In fact, the construction of~$\HM[1]$ generalizes as follows.

Given $n\in\N$, we can consider the unit cube $C^n\subseteq\H^n$. We write $\partial C^n$ for its boundary, which is homotopy equivalent to the $(4n{-}1)$-sphere $\S{4n-1}$. As before, we have an action of~$Q$ on~$C^n$, which restricts to a free action on~$\S{4n-1}$: we have that $\S{4n-1}$ can be expressed as the join of $n$ copies of~$\S3$, and $Q$ acts independently on each of these 3-spheres, as explained in the previous section. We define
\[
  \HM[n]\qdefd\S{4n-1}/Q
\]
as the quotient of the sphere under this action.
By~\cite[Proposition 1.40]{hatcher}, we thus have the fiber sequence
\[
  \begin{tikzcd}
    Q\ar[r]&\S{4n-1}\ar[r]&\HM[n]
  \end{tikzcd}
\]
which simply encodes the fact that the map on the right has $Q$ as fibers, from which we deduce that $\pi_i(\HM[n])=\pi_i(\B Q)$ for $0\leq i\leq 4n-2$, using the induced long exact sequence in homotopy~\cite[Theorem 4.41]{hatcher}.

There is a canonical inclusion $\H^n\to\H^{n+1}$ (obtained by adding $0$ as last coordinate) which induces an inclusion $\partial C^n\to\partial C^{n+1}$, and thus a map $\HM[n]\to\HM[n+1]$. Taking the inductive limit yields a space $\HM[\infty]$ which has the same homotopy groups as $\B Q$ (because homotopy groups commute with inductive limits~\cite[Chapter 9, Section 4]{may}) and is thus homotopy equivalent to $\B Q$~\cite[Theorem 1B.8]{hatcher}. However, this model has the advantage of being a CW-complex which can be explicitly described.

Again, all the constructions above are very ``classical'' in nature, in the sense that we use strict quotients of spaces, boundaries of spaces, and so on. We will see in the rest of the paper how to perform them in a homotopy invariant way, so that they can be formalized in homotopy type theory.
\section{Homotopy type theory}
\label{hott}
The goal of this section is to establish notation and recall the main tools of homotopy type theory used throughout the paper. A comprehensive introduction to the topic can be found in the HoTT Book~\citep{hottbook}.

\subsection{Elementary constructions}
We denote by $\U$ the \emph{universe} of all types (strictly speaking, one should assume a hierarchy of universes, but we will not detail this as it plays no essential role here). Given a type $A:\U$ and a type family~$B:A\to\U$, we write $\Pi A.B$, or $\Pi(x:A).\lapp Bx$, or $(x:A)\to\lapp Bx$ for the associated \emph{dependent product} type: its inhabitants are dependent functions and we write $\labs xt$ for the function that associates to~$x$ (of type~$A$) the term~$t$ (of type~$\lapp Bx$). As usual, we write $A\to B$ for the type of non-dependent functions (here, we suppose that $B$ is a non-dependent type). Similarly, we write $\Sigma A.B$ or $\Sigma(x:A).\lapp Bx$ for the type of \emph{dependent sums}, whose elements are dependent pairs $\pair tu$ consisting of a term~$t$ (of type~$A$) and a term~$u$ (of type $\lapp Bt$). The two dependent projections are respectively written $\fst:\Sigma A.B\to A$ and $\snd:(x:\Sigma A.B)\to\lapp B{(\lapp\fst x)}$. Given two type families $B:A\to\U$ and $B':A'\to\U$, and maps $f:A\to A'$ and $g:(x:A)\to\lapp Bx\to\lapp{B'}{(\lapp fx)}$, we write $\Sigma f.g:\Sigma A.B\to\Sigma A'.B'$ for the canonically induced map.
Given a natural number~$n$, we still write~$n$ for the type with $n$ elements, and even occasionally use the comprehension notation such as $\set{a,b,c}$ to denote the type~$3$ where the three elements are called $a$, $b$ and $c$. In particular, we denote by~$0$ (\resp $1$) the initial (\resp terminal) type.

We write $t\defd u$ to indicate that two terms~$t$ and~$u$ are definitionally equal. Homotopy type theory also provides, for all elements $t$ and $u$ of common type~$A$, a type $t=u$ of \emph{propositional equalities} (or \emph{identities}) between~$t$ and~$u$, also called \emph{paths} because of their topological interpretation. We write $\refl[x]:x=x$ for the trivial path on a point $x:A$. Given paths $p:x=y$ and $q:y=z$, we write $p\pcomp q:x=z$ for the path obtained by concatenation, and $\sym p:y=x$ for the inverse of~$p$.
Given an identity $p:A=B$ between types~$A$ and~$B$, we write $\transport p:A\to B$ for the canonical function induced by transport~\cite[Lemma 2.3.1]{hottbook}.
Given a function $f:A\to B$ and a path $p:t=u$ in~$A$, we write $\ap f(p):f(t)=f(u)$ (or even sometimes $f(p)$) for the canonical path in~$B$~\cite[Lemma 2.2.1]{hottbook} which can be thought of as being obtained by applying~$f$ pointwise to~$p$.
All constructions in homotopy type theory are compatible with identities; in particular, all (co)limits are taken up to homotopy.

A \emph{pointed type} consists of a type~$A$ together with an element of~$A$, which we often denote by~$\pt$, and a pointed map is a function preserving the distinguished element. Given a pointed type~$A$, we write $\Loop A\defd(\pt=\pt)$ for the \emph{loop space} of~$A$.

Given a type~$A$, we write $\ptrunc{A}$ (\resp $\strunc{A}$, \resp $\gtrunc{A}$) for its propositional (\resp set, \resp groupoid) \emph{truncation}, and $\trunq n-:A\to\trunc nA$ for the canonical quotient map. In fact, the $n$-truncation $\trunc nA$ of~$A$ can be defined for any natural number~$n$, and we say that $A$ is \emph{$n$-connected} when $\trunc nA$ is contractible, and we sometimes simply say \emph{connected} for a $0$-connected type.
By extension, a map $f:A\to B$ is \emph{$n$-connected} when the homotopy fiber $\fib fy$ (defined below) is $n$-connected for every $y:B$.

\subsection{Fibrations}
\label{sec:action-fibration}
Any type family~$F:A\to\U$ can be thought of as a fibration with~$A$ as base space, whose fiber at~$x:A$ is simply $\lapp Fx$. The total space of such a fibration is $\Sigma A.F$, which comes equipped with a projection map
\[
  \fst:\Sigma A.F\to A
\]
onto the base space given by the first projection. Conversely, any map $f:X\to A$ (exhibiting~$X$ as a type over~$A$) induces a type family $\fib f:A\to\U$ which associates to each element $y:A$ its \emph{homotopy fiber} defined as
\[
  \fib f y\qdefd\Sigma(x:X).(y=\lapp fx)
\]
which corresponds to the pullback
\[
  \begin{tikzcd}
    \fib fy\ar[d,"i"']\ar[r]\ar[dr,phantom,pos=0,"\lrcorner"]&1\ar[d,"y"]\\
    X\ar[r,"f"']&A
  \end{tikzcd}
\]
These constructions provide an equivalence between the two perspectives on fibrations, as type families indexed by~$A$ or as types over~$A$. This is known as the Grothendieck duality~\cite[Theorem 4.8.3]{hottbook}:

\begin{proposition}[Grothendieck duality]
  \label[proposition]{grothendieck-duality}
  For any type $A:\U$, we have an equivalence between type families and types over~$A$
  \[
    (A \to \U)
    \quad\equivto\quad
    (\Sigma (B:\U).B \to A)
  \]
\end{proposition}

We say that a pair of composable maps
\[
  \begin{tikzcd}
    F\ar[r,"i"]&B\ar[r,"f"]&A
  \end{tikzcd}
\]
with $A$ pointed, is a \emph{fiber sequence} when $F=\fib f\pt$ and $i$ is the canonical map given by the above pullback. This indicates that $f$ can be regarded as a fibration whose fibers are all~$F$ (when $A$ is connected, all the fibers are indeed merely equivalent to $F$). Fiber sequences play a role analogous to short exact sequences in abelian categories.

\label{sec:flattening}

For reasons indicated above, given a type family~$F:A\to\U$, it is often desirable to compute its total space~$\Sigma A.F$. When~$A$ is constructed as a colimit (such as a coequalizer or a pushout), this total space can also be obtained as a similar colimit: this is called the \emph{flattening lemma}. For pushouts, which is the variant we will use here, this fact can be formulated as follows~\cite[Lemma 8.5.3]{hottbook}:

\begin{lemma}[Flattening lemma for pushouts]
  \label[lemma]{flattening-pushout}
  Consider a pushout square
  \[
    \begin{tikzcd}
      X\ar[d,"g"']\ar[r,"f"]\ar[dr,phantom,pos=1,"\ulcorner"]&A\ar[d,"i"]\\
      B\ar[r,"j"']&A\sqcup_XB
    \end{tikzcd}
  \]
  with $p:(x:X)\to(i\circ f)(x)=(j\circ g)(x)$ witnessing its commutativity, and a type family
  \[
    F:A\sqcup_XB\to\U
  \]
  The following square of total spaces is also a pushout:
  \[
    \begin{tikzcd}[column sep=large]
      \Sigma X.(F\circ i\circ f)\ar[d,"\Sigma g.e"']\ar[r,"\Sigma f.(\labs\_{\id{}})"]\ar[dr,phantom,pos=1,"\ulcorner"]&\Sigma A.(F\circ i)\ar[d,"\Sigma i.(\labs\_{\id{}})"]\\
      \Sigma B.(F\circ j)\ar[r,"\Sigma j.(\labs\_{\id{}})"']&\Sigma(A\sqcup_XB).F
    \end{tikzcd}
  \]
  where
  \begin{align*}
    e:(x:X)\to(F\circ i\circ f)(x)&\to(F\circ j\circ g)(x)\\
    x\,y&\mapsto \transport{(\ap F(p\,x))}\,y
  \end{align*}
  is the map canonically induced by~$p$.
\end{lemma}

\noindent
A similar result holds when~$A$ is constructed as a coequalizer~\cite[Lemma 6.12.2]{hottbook}.

\subsection{Groups and actions}
A group $G$ is a set together with a multiplication and unit satisfying the usual axioms. A \emph{delooping} of~$G$ is a pointed connected space $\B G$ together with an isomorphism of groups $\Loop\B G\equivto G$ (where the group structure on $\Loop\B G$ is the one induced by concatenation of paths). A type with this property can be shown to exist and be unique~\citep{warn2023eilenberg}, so that we can talk about \emph{the} delooping of a group. In fact, the loop space provides an equivalence between pointed connected groupoids and groups~\citep{symmetry}. This equivalence implies in particular that we have the following recursion principle:

\begin{proposition}
  \label[proposition]{delooping-rec}
  Given a group~$G$, a pointed groupoid $X$ and a group morphism $f:G\to\Loop X$, there exists a unique pointed map $\tilde f:\B G\to X$ such that the composite
  $
  \begin{tikzcd}[cramped,sep=small]
    G\ar[r]&\Loop\B G\ar[r,"\Loop \tilde f"]& \Loop X
  \end{tikzcd}
  $ is~$f$.
\end{proposition}

A map $\phi:\B G\to\U$ can be interpreted as an \emph{action} of $G$ on the type~$X\defd\phi(\pt)$. Indeed, any $a\in G$ corresponds to a path $a:\pt=\pt$ in~$\B G$, and thus induces, by transport, an endomorphism $\transport{(\ap\phi a)}:X\to X$. In the case where $X$ is a set, this can be shown to correspond to the traditional notion of action. More generally, a map $\B G\to\U$ can be thought of as an action of $G$ on a higher type, in a coherent way. With this point of view, the \emph{homotopy quotient} of~$X$ by~$G$ is the type $X\hq G\eqdef\Sigma(\B G).\phi$. The Grothendieck duality (\cref{grothendieck-duality}), which states that the type of maps $A\to\B G$ is equivalent to the type of families $\B G\to\U$ can thus be read as the following equivalence, which we call the \emph{action-fibration duality}:

\begin{proposition}
  \label[proposition]{action-fibration}
  A fiber sequence
  \[
    \begin{tikzcd}
      X\ar[r]&A\ar[r,"f"]&\B G
    \end{tikzcd}
  \]
  corresponds to a coherent action of~$G$ on~$X$ (\ie a map $\phi:\B G\to\U$ with $\phi(\pt)=X$) with $A$ being the homotopy quotient~$X\hq G$.
\end{proposition}

\noindent
In the above proposition, the action is obtained as the map $\fib f:\B G\to\U$, whose total space is $A=\Sigma(\B G).\fib f$ which is, by definition, the homotopy quotient of~$X$ by the action of~$G$.

\section{The hypercubical manifold in homotopy type theory}
\label{hypercubical-hott}
We now show that the previous constructions can be carried out within the framework of homotopy type theory, which requires more than a straightforward translation. It is clear that we can use HITs to give a direct definition of~$\HM$ as a type consisting of two elements, four identities, three 2-identities and one 3-identity~\citep{dylan}. However, it is not at all obvious that this construction arises from a quotient of~$\S3$ by an action of~$Q$. Defining such an action amounts to constructing a map $\phi:\B Q\to\U$ such that $\phi(\pt{})=\S3$, but the induction principle for $\B Q$ given in \cref{delooping-rec} only allows eliminating into a groupoid, yet $\U$ is not a groupoid, nor is any subtype which contains~$\S3$~\cite[Corollary 6.4.3]{hottbook}.
Instead, since we know that the homotopy quotient of the action encoded by~$\phi$ is $\HM$, the map~$\phi$~corresponds to a fibration $\S3\to\HM\to\B Q$ by the action-fibration duality (\cref{sec:action-fibration}), and this fibration is much easier to define directly. Note that the map $\S3\to\HM$ is the one considered in \cref{hypercubical-manifold}, so that there are eight $n$-cells in the fiber over each $n$-cell in the base. This is encoded by the following long fiber sequence combining the previous ones:
\[
  \begin{tikzcd}
    Q\ar[r]&\S3\ar[r]&\HM\ar[r,"\phi"]&\B Q
  \end{tikzcd}
\]

We will proceed in two steps. We first define a map $\phi:\HM\to\B Q$ using the elimination principle of~$\HM$ (coming from its definition as a HIT). Then, we show that its fiber $\fib\phi(\pt{})$ (which we simply write as $\fib\phi$) is precisely $\S3$ by the flattening lemma. More precisely, writing~$\HM_k$ for the $k$-skeleton of~$\HM$, we reason by induction on~$k$, and compute the fiber of the restriction $\HM_k\to\B Q$ to be the $k$-skeleton of~$\S3$, or more precisely a cubical model of~$\S3$.

\subsection{The hypercubical manifold}
We are now ready to define the type that will play a central role in what follows:

\begin{definition}
  The hypercubical manifold~$\HM$ is the HIT generated by
  \begin{itemize}
  \item 2 elements $a$ and $b$,
  \item 4 identities $x,y,z,w:a=b$,
  \item 3 identities between identities as in \cref{2cells}:
    \begin{align*}
      \alpha:y\pcomp\sym z&=x\pcomp\sym w
      &
      \beta:y\pcomp\sym w&=z\pcomp\sym x
      &
      \gamma:z\pcomp\sym w&=x\pcomp\sym y
    \end{align*}
    which can respectively be pictured as the squares
    \[
      \begin{tikzcd}
        b\ar[dr,phantom,"\alpha"]&\ar[l,"z"']\ar[d,"w"]a\\
        a\ar[u,"y"]\ar[r,"x"']&b
      \end{tikzcd}
      \qquad\qquad\qquad
      \begin{tikzcd}
        b\ar[dr,phantom,"\beta"]&\ar[l,"w"']\ar[d,"x"]a\\
        a\ar[u,"y"]\ar[r,"z"']&b
      \end{tikzcd}
      \qquad\qquad\qquad
      \begin{tikzcd}
        b\ar[dr,phantom,"\gamma"]&\ar[l,"w"']\ar[d,"y"]a\\
        a\ar[u,"z"]\ar[r,"x"']&b
      \end{tikzcd}
    \]
  \item 1 identity between identities between identities as in \cref{hm}.
  \end{itemize}
\end{definition}

We will see in \cref{K-pi1} that the fundamental group of~$\HM$ is~$Q$. This is expected from the topological constructions (see \cref{hypercubical-manifold}), but showing this directly is not easy (although it should be doable using the Seifert-van Kampen theorem~\cite[Section 8.7]{hottbook}).

\begin{lemma}
  \label[lemma]{K-connected}
  The type~$\HM$ is connected.
\end{lemma}
\begin{proof}
  We must show that there merely exists a path from~$a$ to~$x$ for any $x : \HM$. By the induction principle of $\HM$, and the fact that the goal is a proposition, it is sufficient to provide paths of type $a = a$ and $a = b$. We may take $\refl[a]$ and $x$ respectively.
\end{proof}

\subsection{Definition of the action}
We now define the map $\phi:\HM\to\B Q$. To motivate its definition, we can start from the computation of \cref{hm-fundamental-group} which identifies the fundamental group of~$\HM$ with $Q$ (an alternative and more abstract starting point will also be given in \cref{galois-fibration}). We see that the path $x$ corresponds to~$i$, the path $y$ corresponds to~$j$, the relation $z=yx$ imposes that $z$ should correspond to $ji=-k$, and the path $w$ was contracted and should thus correspond to~$1$. This suggests the following definition of the map $\phi$, using the elimination principle of~$\HM$.

\begin{definition}
  \label[definition]{def:phi}
  We define a morphism $\phi:\HM\to\B Q$ on 0-cells by
  \[
    \phi(a)\defd\phi(b)\defd\pt{}
  \]
  and on 1-cells by
  \begin{align*}
    \phi(x)&\defd i
    &
    \phi(y)&\defd j
    &
    \phi(z)&\defd -k
    &
    \phi(w)&\defd 1
  \end{align*}
  The 2-cells are sent to 2-cells coming from the following relations which are known to hold in~$Q$:
  \begin{align*}
    \phi(\alpha)&:i=jk
    &
    \phi(\beta)&:j=ki
    &
    \phi(\gamma)&:k=ij
  \end{align*}
  On 3-cells, $\phi$ is unambiguously defined because $\B Q$ is a groupoid.  
\end{definition}

For $k\in\N$, we define $\HM_k$ to be the $k$-skeleton of $\HM$, \ie $\HM_k$ is defined as the HIT generated by the same $i$-cells for $i\leq k$ (and none for $i>k$). We write $\phi_k:\HM_k\to\B Q$ for the composite map $\begin{tikzcd}[cramped,sep=small]\HM_k\ar[r,hook]&\HM\ar[r,"\phi"]&\B Q\end{tikzcd}$.
Our goal is now to compute the (homotopy) fiber $\fib{\phi_k}$, by induction on~$k$, using the flattening lemma, and ultimately show that $\fib\phi\eqdef\fib{\phi_3}=\S3$ in \cref{the-fiber-sequence}.

\subsection{Fiber of the 0-skeleton}
The 0-skeleton $\HM_0$ is a space consisting of two points~$a$ and~$b$, and the induced map $\phi_0:\HM_0\to\B Q$ is the constant map with $\pt{}$ as image. Its fiber is
\begin{align*}
  \fib{\phi_0}
  &\eqdef
  \Sigma(x:\HM_0).(\pt{}=\phi_0(x))
  \\
  &=
  (\pt{}=\phi(a))
  \sqcup
  (\pt{}=\phi(b))
  \\
  &=
  \set{a,b}\times(\pt{}=\pt{})
  \\
  &=
  \set{a,b}\times Q  
\end{align*}
so that we have a fiber sequence
\[
  \begin{tikzcd}
    \set{a,b}\times Q\ar[r,"\fst"]&\HM_0\ar[r,"\phi_0"]&\B Q
  \end{tikzcd}
\]

\subsection{Fiber of the 1-skeleton}
The 1-skeleton $\HM_1$ is that of the complex pictured in \cref{hm}, with two 0-cells $a$ and $b$, and four 1-cells $x$, $y$, $z$ and~$w$. It can be obtained from $\HM_0$ as the coequalizer
\[
  \begin{tikzcd}
    \set{x,y,z,w}\ar[r,shift left,"\sigma"]\ar[r,shift right,"\tau"']&\HM_0\ar[r,dotted,"\iota_1"]&\HM_1
  \end{tikzcd}
\]
where the maps $\sigma$ and $\tau$ respectively send an edge to its source and target, \ie they are the constant maps respectively equal to~$a$ and $b$.

Our aim is now to compute the fiber of $\phi_1:\HM_1\to\B Q$. We consider the type family
\begin{align*}
  F_1:\HM_1&\to\U\\
  x&\mapsto(\pt{}=\phi_1(x))
\end{align*}
This family is defined so that we have
\[
  \fib{\phi_1}\eqdef\Sigma(x:\HM_1).F_1(x)
\]
which can be computed using the flattening lemma for coequalizers (see \cref{sec:flattening} and \cite[Section 6.12]{hottbook}) to be the coequalizer
\begin{equation}
  \label{hm1-coeq}
  \begin{tikzcd}[sep=large]
    \Sigma\set{x,y,z,w}.(F_1\circ\iota_1\circ\sigma)
    \ar[r,shift left,"\Sigma\sigma.\lambda\_\id{}"]
    \ar[r,shift right,"\Sigma\tau.e"']
    &\Sigma\HM_0.(F_1\circ\iota_1)
    \ar[r,dotted,"\Sigma\iota_1.\lambda\_\id{}"]
    &\Sigma\HM_1.F_1
  \end{tikzcd}
\end{equation}
where $e$ is defined as
\begin{align*}
  e:(p:\set{x,y,z,w})\to(F_1\circ\iota_1\circ\sigma)(p)&\to(F_1\circ\iota_1\circ\tau)(p)\\
  p\ q&\mapsto\transport{(\ap{F_1}\,p)}\,q
\end{align*}
(note that, above, both types $(F_1\circ\iota_1\circ\sigma)(p)$ and $(F_1\circ\iota_1\circ\tau)(p)$ are equal to $\Loop\B Q$).
We have $F_1(\iota_1(\sigma(x)))\eqdef(\pt{}=\pt{})=Q$ and similarly for $y$, $z$ and~$w$, so that \eqref{hm1-coeq} can be rewritten as
\begin{equation}
  \label{hm1-coeq'}
  \begin{tikzcd}[sep=large]
    \set{x,y,z,w}\times Q
    \ar[r,shift left,"\Sigma\sigma.\lambda\_\id{}"]
    \ar[r,shift right,"\Sigma\tau.e"']
    &\HM_0\times Q
    \ar[r,dotted,"\Sigma\iota_1.\lambda\_\id{}"]
    &\fib{\phi_1}
  \end{tikzcd}
\end{equation}
By properties of transport in path types~\cite[Theorem 2.11.4]{hottbook}, we thus have $e\,p\,q=q\pcomp\ap\phi\,p$. Equivalently, in \eqref{hm1-coeq'}, we can thus define
\begin{align*}
  e:\set{x,y,z,w}&\to Q\to Q\\
  p\ q&\mapsto q\cdot\ap\phi(p)
\end{align*}
by slightly overloading~$e$, where $\cdot$ denotes the multiplication in~$Q$, and we implicitly identify $\ap\phi(p)$ which is an element of~$\Loop\B Q$ with the corresponding element of~$Q$.

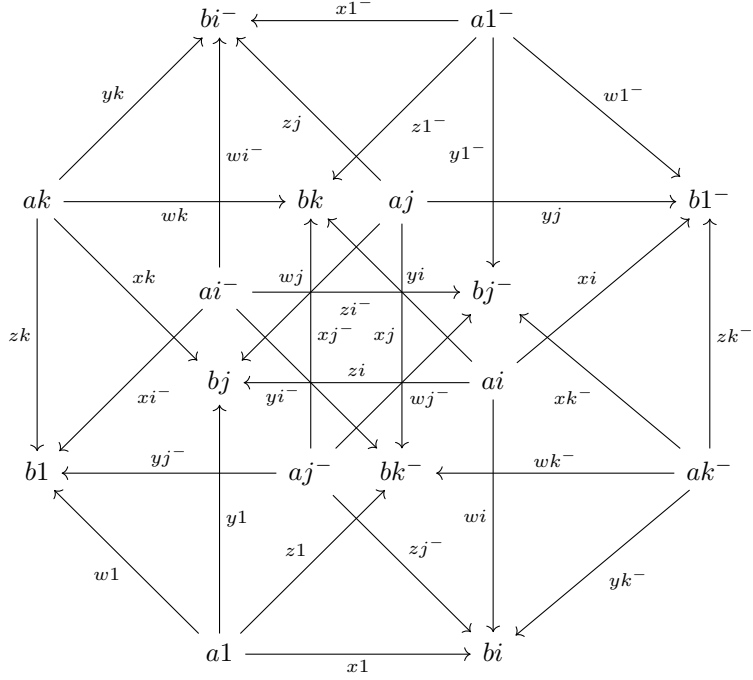
\begin{figure}[t]
  \centering
  \[
  \begin{tikzcd}[column sep={12mm,between origins},row sep={12mm,between origins}]
    &&bi^-&&&a1^-\ar[lll,"x1^-"']\ar[ddll,"z1^-"]\ar[ddd,"y1^-"']\ar[ddrr,"w1^-"]\\
    \\
    ak\ar[ddd,"zk"']\ar[uurr,"yk"]\ar[rrr,"wk"']\ar[ddrr,"xk"]&&&bk&aj\ar[uull,"zj"]\ar[ddll,"wj"']\ar[ddd,"xj"']\ar[rrr,"yj"']&&&b1^-\\
    &&ai^-\ar[ddll,"xi^-"]\ar[uuu,"wi^-"']\ar[rrr,"zi^-"']\ar[ddrr,"yi^-"']&&&bj^-\\
    &&bj&&&ai\ar[uull,"yi"']\ar[lll,"zi"']\ar[ddd,"wi"']\ar[uurr,"xi"]\\
    b1&&&aj^-\ar[lll,"yj^-"']\ar[uuu,"xj^-"']\ar[uurr,"wj^-"']\ar[ddrr,"zj^-"]&bk^-&&&ak^-\ar[uull,"xk^-"]\ar[lll,"wk^-"']\ar[ddll,"yk^-"]\ar[uuu,"zk^-"']\\
    \\
    &&\ar[uull,"w1"]a1\ar[uurr,"z1"]\ar[uuu,"y1"']\ar[rrr,"x1"']&&&bi
  \end{tikzcd}
\]
  \caption{The fiber of~$\phi_1$ (the 1-skeleton of a cubical model of~$\S3$).}
  \label[figure]{tesseract}
\end{figure}

The coequalizer \eqref{hm1-coeq'} can be interpreted as a description of $\fib{\phi_1}$ as a 1\nbd-dimen\-sional cell complex, with $\HM_0\eqdef\set{a,b}\times Q$ as 0-cells and $\set{x,y,z,w}\times Q$ as 1\nbd-cells, the source and target of a 1-cell $(p,q):\set{x,y,z,w}\times Q$ being respectively $(a,q)$ and $(b,q\cdot\ap\phi(p))$. This complex is pictured in \cref{tesseract}. Here, we write $ai$ (\resp $ai^-$) for $(a,i)$ (\resp $(a,-i)$) and similarly for other cells, so that the general form of a 1-cell is
\[
  \begin{tikzcd}
    aq\ar[r,"pq"]&b(q\cdot\ap\phi(p))
  \end{tikzcd}
\]
We can observe that this is the 1-skeleton of a tesseract (a 4-dimensional cube), which is homotopic to $\S3$ as expected.

\begin{remark}
  \begin{figure}[t]
    \centering
    \[
      \begin{tikzcd}
        ai\ar[rrr,"xi"]\ar[ddr,"yi"']&&&\ar[dll,"y1^-"']a1^-\ar[ddd,"x1^-"]\\
        &aj^-\ar[ddl,"yj^-"']\ar[d,"xj^-"]&\ar[l,"xk^-"]\ar[ull,"yk^-"']ak^-\\
        &ak\ar[drr,"yk"']\ar[r,"xk"]&\ar[u,"xj"]aj\ar[uur,"yj"']\\
        a1\ar[urr,"y1"']\ar[uuu,"x1"]&&&\ar[lll,"xi^-"]\ar[uul,"yi^-"']ai^-
      \end{tikzcd}
    \]
    \caption{Cayley graph of~$Q$ generated by $i$ and $j$.}
    \label[figure]{cayley}
  \end{figure}
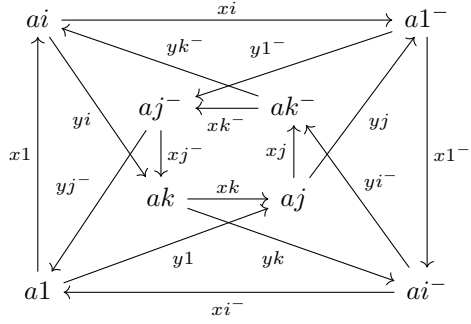
  The group~$Q$ is generated by~$i$ and~$j$, and so it admits the Cayley graph~$C$ shown in \cref{cayley}. This graph describes a HIT (with the vertices as generators and edges as identities) which can be obtained as the fiber of the canonical map $\B\set{i,j}^*\to\B Q$, where $\set{i,j}^*$ is the free group on two generators, see~\cite[Theorem 16]{delooping-generated}. Similarly, we have a Cayley graph~$C'$ associated to~$\set{x,y,z}$ as generators (where $x$, $y$ and $z$ respectively correspond to $i$, $j$ and $-k$).
  We have a map $f:\HM_1\to\B\set{x,y,z}^*$ which expresses that the target can be obtained from the source by contracting the 1-cell~$w$ to a point, \ie
  \begin{align*}
    f(a)&\eqdef\pt
    &
    f(b)&\eqdef\pt
    &
    f(x)&\eqdef x
    &
    f(y)&\eqdef y
    &
    f(z)&\eqdef z
    &
    f(w)&\eqdef\refl[\pt]
  \end{align*}
  We also have a map $g:\B\set{i,j}^*\to\B\set{x,y,z}^*$ which is the canonical inclusion of cells, respectively sending~$x$ and~$y$ to $i$ and~$j$. These maps make the two triangles on the right commute:
  \[
    \begin{tikzcd}
      \fib{\phi_1}\ar[d,dotted,->>,"f'"']\ar[r]&\HM_1\ar[d,->>,"f"']\ar[dr,"\phi_1"]\\
      C'\ar[r]&\B\set{x,y,z}^*\ar[r]&\B Q\\
      C\ar[u,hook,dotted,"g'"]\ar[r]&\B\set{i,j}^*\ar[u,hook,"g"]\ar[ur]
    \end{tikzcd}
  \]
  where $\twoheadrightarrow$ (\resp $\hookrightarrow$) indicates a map which is surjective (\resp injective) on loop spaces.
  In other words, these maps express that $\B\set{i,j}^*$ can be obtained as a ``subquotient'' of $\HM_1$ (quotienting by~$w$ and forgetting $z$). By the universal property of~$C'$ which is a fiber, we have similar maps $f'$ and $g'$ as shown above between the fibers, expressing that~$C$ can be obtained as a subquotient of~$\fib{\phi_1}$, by quotienting all 1-cells of the form $(w,q)$ with $q\in Q$ (and thus identifying all 0-cells $(b,q)$ with $(a,q)$) and forgetting all the 1-cells of the form $(z,q)$ with~$q\in Q$.
  It can be checked that if we perform these operations on~$\HM_1$ pictured in~\cref{tesseract} we recover the traditional Cayley graph pictured in~\cref{cayley}.
\end{remark}


\subsection{Fiber of the 2-skeleton}
We write $\square$ for the HIT corresponding to an empty square, generated by four 0-cells $v_0,v_1,v_2,v_3$ and four 1-cells $e_0,e_1,e_2,e_3$ attached as in the following picture:
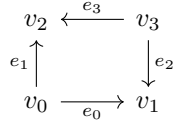
\[
  \begin{tikzcd}
    v_2&\ar[l,"e_3"']v_3\ar[d,"e_2"]\\
    v_0\ar[u,"e_1"]\ar[r,"e_0"']&v_1
  \end{tikzcd}
\]
We also write $\blacksquare$ for the HIT corresponding to a filled square, obtained from the previous one by attaching a 2-cell as expected. There is an induced canonical inclusion $\square\into\blacksquare$. The type $\blacksquare$ is contractible, and we will consider a particular equivalence $\blacksquare\equivto 1$ consisting of the terminal map $\tau:\blacksquare\to 1$ and the map $\upsilon:1\to\blacksquare$ pointing at $v_0$. Since these form an equivalence, their composite is the identity on~$\blacksquare$, which is witnessed by a homotopy $\varepsilon:(x:\blacksquare)\to v_0=x$. Up to deformation (being an equivalence is a property), we can consider a particular definition for $\varepsilon$, which is obtained by induction on $\blacksquare$, and satisfies
\begin{align*}
  \varepsilon(v_0)&\defd\refl
  &
  \varepsilon(v_1)&\defd e_0
  &
  \varepsilon(v_2)&\defd e_1
  &
  \varepsilon(v_3)&\defd e_0\pcomp\overline{e_2}
\end{align*}
definitionally (on 1- and 2-cells, this is defined arbitrarily, using the fact that the filled square is contractible). This will slightly simplify subsequent computations.

The 2-skeleton $\HM_2$ can be obtained from $\HM_1$ as the pushout on the left
\[
  \begin{tikzcd}
    \set{\alpha,\beta,\gamma}\times\square\ar[d,"\psi"']\ar[r,hook,"\iota"]\ar[dr,phantom,pos=1,"\ulcorner"]&\set{\alpha,\beta,\gamma}\times\blacksquare\ar[d,"\tilde\rho"]\\
    \HM_1\ar[r,hook,"\iota_2"']\ar[ur,equals,"\tilde h",shorten=25pt]&\HM_2
  \end{tikzcd}
  \qquad\qquad\quad
  \begin{tikzcd}
    \set{\alpha,\beta,\gamma}\times\square\ar[d,"\psi"']\ar[r,hook,"\fst"]\ar[dr,phantom,pos=1,"\ulcorner"]&\set{\alpha,\beta,\gamma}\ar[d,"\rho"]\\
    \HM_1\ar[r,hook,"\iota_2"']\ar[ur,equals,"h",shorten=20pt]&\HM_2
  \end{tikzcd}
\]
where the vertical map $\psi$ sends each formal cell to its boundary, as pictured in \cref{2cells}, and the horizontal map $\iota$ is the identity on the first component and the canonical inclusion of the boundary of the square. We write $\tilde h:\tilde\rho\circ\iota=\iota_2\circ\psi$ for the canonical homotopy witnessing the commutativity of the square.
By composing the upper-right corner with the equivalence $\blacksquare\equivto 1$ described above, we obtain the pushout figured on the right above, where $\fst:\set{\alpha,\beta,\gamma}\times\square\to\set{\alpha,\beta,\gamma}$ is the first projection $\fst\eqdef(\id{\set{\alpha,\beta,\gamma}}\times\tau)\circ\iota$ and $\rho:\set{\alpha,\beta,\gamma}\to\HM_2$, defined as $\rho\eqdef\tilde\rho\circ(\id{\set{\alpha,\beta,\gamma}}\times\upsilon)$ sends $\alpha$ to the corner of the square $\alpha$ in~$\HM_2$, corresponding to $v_0$, \ie $a$ (and similarly for $\beta$ and $\gamma$). Hence $\rho$ is the constant map equal to $a$. More interestingly, the family of equalities $\tilde h$ is transformed into the family $h:\rho\circ\fst=\iota_2\circ\psi$ obtained as the composite $h=\tilde\varepsilon\pcomp\tilde h$ where $\tilde\varepsilon:\rho\circ\fst\eqdef\tilde\rho\circ(\id{\set{\alpha,\beta,\gamma}}\times(\upsilon\circ\tau))\circ\iota=\tilde\rho\circ\iota$ is the equality canonically induced by $\varepsilon$.
%

As previously, we consider the type family $F_2:\HM_2\to\U$ defined by
\[
  F_2(x)\qdefd(\pt{}=\phi_2(x))
\]
The flattening lemma for pushouts (\cref{flattening-pushout}) ensures that $\Sigma\HM_2.F_2$ can be computed as the pushout
\begin{equation}
  \label{pushout2}
  \begin{tikzcd}[column sep=large]
    \Sigma(\set{\alpha,\beta,\gamma}\times\square).(F_2\circ\rho\circ\fst)\ar[d,"\Sigma\psi.E"']\ar[r,"\Sigma\fst.\id{}"]\ar[dr,phantom,pos=1,"\ulcorner"]&\Sigma\set{\alpha,\beta,\gamma}.(F_2\circ\rho)\ar[d,"\Sigma\rho.\lambda\_.\id{}"]\\
    \Sigma\HM_1.(F_2\circ\iota_2)\ar[r,"\Sigma\iota_2.\lambda\_\id{}"']&\Sigma\HM_2.F_2
  \end{tikzcd}
\end{equation}
The morphism
\[
  E:((\omega,s):\set{\alpha,\beta,\gamma}\times\square)\to(F_2\circ\rho\circ\fst)(\omega,s)\to(F_2\circ\iota_2\circ\psi)(\omega,s)
\]
is defined as follows on $(\omega,s):\set{\alpha,\beta,\gamma}\times\square$. We have a path
\[
  h(\omega,s):(\rho\circ\fst)(\omega,s)=(\iota_2\circ\psi)(\omega,s)
\]
which, by applying~$F_2$, induces a path
\[
  \ap F_2(h(\omega,s)):(F_2\circ\rho\circ\fst)(\omega,s)=(F_2\circ\iota_2\circ\psi)(\omega,s)
\]
which, by transport, induces the desired map
\[
  \transport{(\ap F_2(h(\omega,s)))}:(F_2\circ\rho\circ\fst)(\omega,s)\to(F_2\circ\iota_2\circ\psi)(\omega,s)
\]
Now, observe that for any $\omega:\set{\alpha,\beta,\gamma}$, we have
\[
  (F_2\circ\rho)(\omega)
  \qeq
  F_2(a)
  \qeq
  (\pt = \phi_2(a))
  \qeq
  (\pt =\pt)
  \qeq
  Q 
\]
because $\rho$ is the constant map at $a : \HM_2$, and thus the upper-right corner of the pushout \cref{pushout2} is
\[
 \Sigma(\omega:\set{\alpha,\beta,\gamma}).(F_2\circ\rho)(\omega) = \set{\alpha,\beta,\gamma} \times Q
\]
the upper-left corner can be transformed similarly as $\Sigma(\set{\alpha,\beta,\gamma}\times\square).Q$, which is equivalent to $\set{\alpha,\beta,\gamma}\times Q\times\square$. Finally, the lower left corner
\[
  \Sigma\HM_1.(F_2\circ\iota_2)\qeqdef\Sigma\HM_1.F_1
\]
is precisely the fiber $\fib{\phi_1}$ computed in the previous section. The pushout \cref{pushout2} can thus be slightly simplified as the pushout
\[
  \begin{tikzcd}
    \set{\alpha,\beta,\gamma}\times Q\times\square\ar[d,"\tilde\psi"']\ar[r,"\fst"]\ar[dr,phantom,pos=1,"\ulcorner"]&\set{\alpha,\beta,\gamma}\times Q\ar[d,"\Sigma\rho.\lambda\_.\id{}"]\\
    \fib{\phi_1}\ar[r,"\Sigma\iota_2.\lambda\_\id{}"']&\Sigma\HM_2.F_2
  \end{tikzcd}
\]
The map $\tilde\psi$ is defined, for $\omega:\set{\alpha,\beta,\gamma}$, $q:Q$ and $s:\square$, by
\[
  \tilde\psi(\omega,q,s)
  \quad=\quad
  (\psi(\omega,s),\transport{(\ap F_2(e(\omega,s)))}\,q)
  \quad=\quad
  (\psi(\omega,s),q\cdot\ap\phi_2(e(\omega,s)))
\]
The values of this map can be computed explicitly. We fix $\omega\defd\alpha$ (computing the images for $\beta$ and $\gamma$ is similar), and compute the image depending on~$s$. When $s$ is a vertex, we have:
\begin{itemize}
\item for $v_0$, the image is $(a,q)$,
\item for $v_1$, the image is $(b,q\cdot\phi_2(x))\eqdef(b,q\cdot i)$ because $x$ is the path in~$\alpha$ corresponding to $e_0\eqdef\varepsilon(v_1)$,
\item for $v_2$, the image is $(b,q\cdot\phi_2(y))\eqdef(b,q\cdot j)$ because $y$ is the path in~$\alpha$ corresponding to $e_1\eqdef\varepsilon(v_2)$,
\item for $v_3$, the image is $(a,q\cdot\phi_2(x)\pcomp\phi_2(\ol w))\eqdef(a,q\cdot i)$ because $x\pcomp\ol w$ is the path in~$\alpha$ corresponding to $e_0\pcomp\overline{e_2}\eqdef\varepsilon(v_3)$.
\end{itemize}
The image of an edge $e_i:v_j\to v_k$ is the edge whose first component is the edge in~$\alpha$ corresponding to $e_i$ and the second component is the element of $Q$ corresponding to the first component of the image of its source $v_j$.
The squares $\alpha$, $\beta$ and $\gamma$ are thus respectively sent to the following squares in $\fib{\phi_1}$:
\[
  \begin{tikzcd}
    (b,q\cdot j)&\ar[l,"{(z,q\cdot i)}"']\ar[d,"{(w,q\cdot i)}"](a,q\cdot i)\\
    (a,q)\ar[u,"{(y,q)}"]\ar[r,"{(x,q)}"']&(b,q\cdot i)
  \end{tikzcd}
  \qquad\qquad
  \begin{tikzcd}
    (b,q\cdot j)&\ar[l,"{(w,q\cdot j)}"']\ar[d,"{(x,q\cdot j)}"](a,q\cdot k^-\cdot i^-)\\
    (a,q)\ar[u,"{(y,q)}"]\ar[r,"{(z,q)}"']&(b,q\cdot k^-)
  \end{tikzcd}
\]
and
\[
  \begin{tikzcd}
    (b,q\cdot k^-)&\ar[l,"{(w,q\cdot k^-)}"']\ar[d,"{(y,q\cdot k^-)}"](a,q\cdot i\cdot j^-)\\
    (a,q)\ar[u,"{(z,q)}"]\ar[r,"{(x,q)}"']&(b,q\cdot i)
  \end{tikzcd}
\]
We respectively write $(\alpha,q)$, $(\beta,q)$ and $(\gamma,q)$ for these cells.
Finally, the previous pushout states that $\fib{\phi_2}$ can be obtained from $\fib{\phi_1}$ by attaching cells along these boundaries, for any $q:Q$.

\subsection{Fiber of the 3-skeleton}
The computation of the 3-skeleton can be performed similarly, as we now explain (with fewer details, the computations being similar). We write $\cube$ for the standard cube which is the HIT pictured on the left
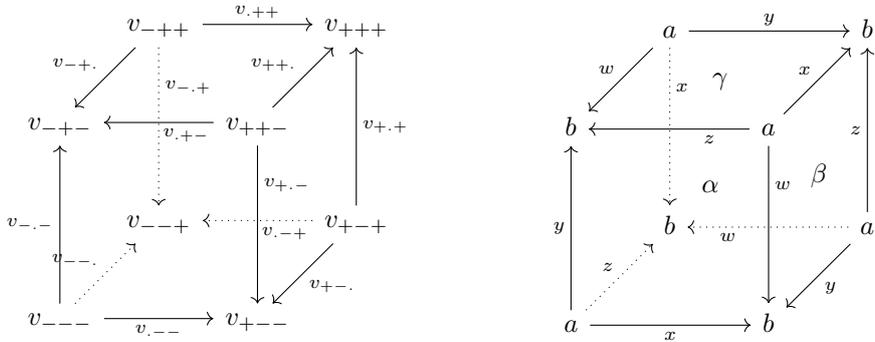
\[
  \begin{tikzcd}[column sep={13mm,between origins},row sep={13mm,between origins}]
    &v_{-++}\ar[dl,"v_{-+.}"']\ar[dd,dotted,"v_{-.+}",near start]\ar[rr,"v_{.++}"]&&v_{+++}\\
    v_{-+-}&&v_{++-}\ar[ll,"v_{.+-}",near start]\ar[dd,"v_{+.-}",near start]\ar[ur,"v_{++.}"]\\
    &v_{--+}&&v_{+-+}\ar[dl,"v_{+-.}"]\ar[ll,dotted,"v_{.-+}",near start]\ar[uu,"v_{+.+}"']\\
    v_{---}\ar[uu,"v_{-.-}"]\ar[ur,dotted,"v_{--.}",pos=.75]\ar[rr,"v_{.--}"']&&v_{+--}
  \end{tikzcd}
  \qquad\qquad\quad
  \begin{tikzcd}[column sep={13mm,between origins},row sep={13mm,between origins}]
    &\ar[dl,"w"']a\ar[dr,phantom,"\gamma"]\ar[dd,dotted,"x",near start]\ar[rr,"y"]&&b\\
    b&&\ar[ll,"z",near start]a\ar[dd,"w",near start]\ar[ur,"x"]\\
    &b&&\ar[ll,dotted,"w",near end]\ar[dl,"y"]a\ar[uu,"z"]\\
    a\ar[uu,"y"]\ar[ur,dotted,"z"]\ar[rr,"x"']\ar[uurr,phantom,"\alpha",near end]&&b\ar[uuur,phantom,"\beta"]
  \end{tikzcd}
\]
we also write $\filledcube$ for the variant with a 3-cell filling the interior of the cube.
We have a map $\chi:\cube\to\HM_2$ sending the standard cube (on the left) to the expected cube in~$\HM_2$ (on the right above), where $\alpha$, $\beta$ and $\gamma$ are respectively the front, right and top faces of the cube. The type~$\HM$ can be obtained as the pushout on the left:
\[
  \begin{tikzcd}
    \cube\ar[d,"\chi"']\ar[r]\ar[dr,phantom,"\ulcorner",pos=1]&\filledcube\ar[d]\\
    \HM_2\ar[r,"\iota_3"']&\HM
  \end{tikzcd}
  \qquad\qquad\qquad\qquad
  \begin{tikzcd}
    \cube\ar[d,"\chi"']\ar[r,"\tau"]\ar[dr,phantom,"\ulcorner",pos=1]&1\ar[d,"\rho"]\\
    \HM_2\ar[r,"\iota_3"']&\HM
  \end{tikzcd}
\]
Since the type $\filledcube$ is contractible, the pushout is equivalent to the one on the right, where the right vertical map is pointing to~$a$ (the image of $v_{---}$).

Consider the type family $F:\HM\to\U$ defined by $F(x)\defd(\pt{}=\phi(x))$. By the flattening lemma for pushouts, we have a pushout as on the left, which can equivalently be computed as on the right:
\[
  \begin{tikzcd}
    \Sigma\cube.(F\circ\rho\circ\tau)\ar[d]\ar[r]\ar[dr,phantom,"\ulcorner",pos=1]&\Sigma 1.(F\circ\rho)\ar[d]\\
    \Sigma \HM_2.(F\circ\iota_3)\ar[r]&\Sigma \HM.F    
  \end{tikzcd}
  \qquad\qquad\qquad
  \begin{tikzcd}
    \cube\times Q\ar[d,"\tilde\chi"']\ar[r]\ar[dr,phantom,"\ulcorner",pos=1]&Q\ar[d]\\
    \fib{\phi_2}\ar[r]&\fib\phi
  \end{tikzcd}
\]
thus, we have an explicit construction of $\fib\phi$.
The map $\tilde\chi$ can be computed in a similar way to the above and sends, for $q:Q$, the cube $\cube$ to
\[
  \begin{tikzcd}
    &\ar[dl,"{(w,q\cdot j)}"'](a,q\cdot j)\ar[dd,dotted,"{(x,q\cdot j)}",near start]\ar[rr,"{(y,q\cdot j)}"]&&(b,q\cdot 1^-)\\
    (b,q\cdot j)\ar[urrr,phantom,"{(\gamma,q\cdot i)}",pos=.6]&&\ar[ll,"{(z,q\cdot i)}",near start](a,q\cdot i)\ar[dd,"{(w,q\cdot i)}",near start]\ar[ur,"{(x,q\cdot i)}"]\\
    &(b,q\cdot k^-)&&\ar[ll,dotted,"{(w,q\cdot k^-)}",near end]\ar[dl,"{(y,q\cdot k^-)}"](a,q\cdot k^-)\ar[uu,"{(z,q\cdot k^-)}"']\\
    (a,q)\ar[uu,"{(y,q)}"]\ar[ur,dotted,"{(z,q)}"]\ar[rr,"{(x,q)}"']\ar[uurr,phantom,"{(\alpha,q)}",near end]&&(b,q\cdot i)\ar[uuur,phantom,"{(\beta,q\cdot k^-)}",pos=.6]
  \end{tikzcd}
\]

Finally, this provides us with a description of $\fib\phi$ as a HIT consisting of cubical cells (16 0-cells, 32 1-cells, 24 2-cells, 8 3-cells) which can easily be checked to form an empty 4-cube, which is homotopy equivalent to $\S3$. We have also provided a computer-checked formalization of this result, consisting of an explicit construction of the equivalence between the above HIT and the standard tesseract in the proof assistant Agda~\citep{git}. In the future, it would be interesting to perform a complete formalization of the above computations, although this is expected to be a considerable amount of work.

\subsection{The fundamental fiber sequence}
The previous computation can be summarized by the following fiber sequence.

\begin{theorem}
  \label[theorem]{the-fiber-sequence}
  We have a fiber sequence
  \[
    \begin{tikzcd}
      \S3\ar[r]&\HM\ar[r,"\phi"]&\B Q
    \end{tikzcd}
  \]
  that we call the \emph{fundamental fiber sequence}, which encodes a coherent action of~$Q$ on~$\S3$ whose homotopy quotient is $\HM$.
\end{theorem}

\noindent
This justifies considering our first definition of $\HM$ as a HIT to be an appropriate definition of the hypercubical manifold (up to homotopy).

As a consequence, we have the following:

\begin{lemma}
  \label[lemma]{phi-2-connected}
  The map $\phi$ is 2-connected.
\end{lemma}
\begin{proof}
  By \cref{the-fiber-sequence}, we have $\trunc2{\fib\phi}=\trunc2{\S3}=1$ since $\S3$ is 2-connected~\cite[Corollary 8.2.2]{hottbook}.
\end{proof}


\begin{proposition}
  \label[proposition]{K-pi1}
  The fundamental group of the hypercubical manifold is $\pi_1(\HM)=Q$.
\end{proposition}
\begin{proof}
  The map $\phi$ is 2-connected by \cref{phi-2-connected}; it thus induces an isomorphism on fundamental groups~\cite[Corollary 8.8.5]{hottbook}, and the fundamental group of~$\B Q$ is $Q$.
\end{proof}


\subsection{Recovering the Galois fibration}
\label{galois-fibration}
Taking groupoid truncations induces a commuting square
\[
  \begin{tikzcd}
    \HM\ar[d,"\gtrunq-"']\ar[r,"\phi"]&\B Q\ar[d,"\gtrunq-"]\\
    \gtrunc{\HM}\ar[r,"\gtrunc{\phi}"']&\gtrunc{\B Q}
  \end{tikzcd}
\]
where the vertical map on the right is an equivalence because, by definition, $\B Q$ is a groupoid. Moreover, since $\phi$ is 1-connected (by \cref{phi-2-connected}), the horizontal map at the bottom is an equivalence~\cite[Lemma 7.5.14]{hottbook}. Thus, up to an automorphism of the codomain, the map $\phi:\HM\to\B Q$ coincides with the map $\gtrunq{{-}}:\HM\to\gtrunc{\HM}$, which is known as the \emph{Galois fibration}~\citep{galoisfib}:

\begin{lemma}
  \label{phi-galois}
  We have a commuting triangle
  \[
    \begin{tikzcd}
      &\ar[dl,"\gtrunq{{-}}"']\HM\ar[dr,"\phi"]\\
      \gtrunc{\HM}\ar[rr,"\sim"']&&\B Q
    \end{tikzcd}
  \]
  where the bottom map is an equivalence.
\end{lemma}

Let us briefly explain the importance of this fibration. Given a pointed connected type~$A$, we have that $\gtrunc{A}$ is a connected groupoid whose fundamental group is $\pi_1(A)$ since we have
\[
  \Trunc0{\gtrunq\pt=_{\gtrunc{A}}\gtrunq\pt}
  \qeq
  \Trunc0{\Trunc0{\pt=_A\pt}}
  \qeq
  \Trunc0{\pt=_A\pt}
\]
where the first equality is~\cite[Theorem 7.3.12]{hottbook} and the second one is~\cite[Corollary 7.3.7]{hottbook}. The fiber of the truncation map $\gtrunq{{-}}:A\to\gtrunc{A}$ is the type on the left below
\[
  \Sigma(x:A).(\gtrunq{\pt{}}=\gtrunq{\,x\,})
  \qeq
  \Sigma(x:A).\strunc{\pt{}=x\,}
\]
which, by \cite[Theorem 7.3.12]{hottbook} again, is equal to the type on the right: this latter type is precisely the universal cover of~$A$ along with the first projection $\fst:\tilde A\to A$ as universal covering map. In other words, we have a fiber sequence
\[
  \begin{tikzcd}
    \tilde A\ar[r,"\fst"]&A\ar[r,"\gtrunq-"]&\B\pi_1(A)
  \end{tikzcd}
\]
which, under the action-fibration duality (\cref{action-fibration}), corresponds to the canonical action of~$\pi_1(A)$ on the universal cover.
The fiber sequence of \cref{the-fiber-sequence} is essentially of this form, with $A\defd\HM$ so we have shown:

\begin{theorem}
  The universal cover of the hypercubical manifold is $\S3$.
\end{theorem}

\noindent
Note that we have defined the map $\phi:\HM\to\B Q$ in \cref{def:phi} by making an educated guess, but we could actually have computed it from the map $\gtrunq{{-}}:\HM\to\gtrunc{\HM}$.

\section{Higher dimensional hypercubical manifolds}
\label{higher}
We now construct higher dimensional variants of the hypercubical manifold: for every natural number $n$, we define a type $\HM[n]$ equipped with a canonical map $\phi^n:\HM[n]\to\B Q$. In particular, for~$n=1$, we recover the previous constructions: $\HM[1]\eqdef\HM$ and $\phi^1\eqdef\phi$. These will be such that there is a fiber sequence
\[
  \begin{tikzcd}
    \S{4n-1}\ar[r]&\HM[n]\ar[r,"\phi^n"]&\B Q
  \end{tikzcd}
\]
exhibiting the fact that $\HM[n]$ is a quotient of $\S{4n-1}$ by an action of $Q$, and therefore the map $\phi^n$ will be $(4n{-}2)$-connected. In this sense, the types $\HM[n]$ provide better and better approximations of $\B Q$ as $n$ increases.

In order to define these in type theory, we should draw inspiration from the topological definitions recalled in \cref{higher-topological-hypercubical}. For instance, $\HM[2]$ is obtained as a quotient of~$\S7$ by~$Q$. Recalling that~$\S7$ can be constructed as a join $\S7=\S3\join\S3$ of two copies of~$\S3$, the action of~$Q$ on~$\S7$ can also be decomposed as a form of join of the action of $Q$ on~$\S3$ used to define~$\HM$. This suggests defining~$\phi^2$ as a join $\phi^2\defd\phi\join\phi$, and similarly for $\phi^n$ with arbitrary~$n$ (we will recall below how the join operation translates in type theory).
This construction is reminiscent of the one of real projective spaces~\citep{buchholtz2017real} and lens spaces~\citep{lens}, which are also obtained as joins of suitable maps to the delooping of their fundamental group.

\subsection{The join operation in type theory}
%
Given two types $A$ and $B$, their \emph{join} $A\join B$, see~\cite[Section~6.8]{hottbook}, is the pushout of the product projections:
\[
  \begin{tikzcd}
    A\times B\ar[d,"\fst"']\ar[r,"\snd"]\ar[dr,phantom,pos=1,"\ulcorner"]&B\ar[d,dotted,"\sndinj"]\\
    A\ar[r,dotted,"\fstinj"']&A\join B
  \end{tikzcd}
\]
In type theory, the type $A\join B$ can be described as the HIT with $A \to A\join B$ and $B \to A\join B$ as base constructors and $(a,b):A\times B\to A\join B$ as a higher constructor adding paths $a=b$ indexed by pairs $(a,b):A\times B$. 

The operation $\S0\join-$ is the suspension operation so that $\S0 \join \S n = \S{n+1}$. The join operation is associative, commutative and admits $0$ as neutral element~\cite[Section~1.8]{brunerie2016homotopy}, and we thus have the fact that $\S n \join \S m = \S{n+m+1}$. It can be shown that joins increase connectivity in the following sense~\cite[Proposition 3]{brunerie2019james}:

\begin{proposition}
  \label[proposition]{join-connected}
  If~$A$ is $m$-connected and~$B$ is $n$-connected, $A\join B$ is $(m{+}n{+}2)$-connected.
\end{proposition}

We can generalize the join operation to pairs of morphisms with common codomain as follows. Given maps~$f:A\to C$ and $g:B\to C$, their \emph{join} $f\join g$ is the universal map obtained from the pushout of the dependent projections~$\fst$ and~$\snd$ from their pullback $A\times_CB$:
\[
  \begin{tikzcd}
    A\times_CB\ar[d,"\fst"']\ar[r,"\snd"]&B\ar[d,"\sndinj"]\ar[ddr,bend left,"g"]\\
    A\ar[drr,"f"',bend right=20]\ar[r,"\fstinj"']&A\join_CB\ar[dr,"f\join g"description]\\
    &&C
  \end{tikzcd}
\]
The source $A\join_CB$ can be thought of as a ``dependent join over~$C$'' in the sense that we have $A\join_1B=A\join B$: we recover the join operation on types by considering terminal maps. An important result about the join operation is the following~\cite[Theorem 2.3.15]{rijke2018classifying}:

\begin{proposition}
  \label[proposition]{fiber-join}
  Given maps $f:A\to C$ and $g:B\to C$, for any $x:C$ we have an equality
  \[
    \fib{f\join g} x = (\fib f x) \join (\fib g x)
  \]
\end{proposition}

\noindent
For this reason $A\join_CB$ is also sometimes called the \emph{fiberwise join} of~$A$ and~$B$ over~$C$.
As a direct consequence of the previous proposition, we have the following result:

\begin{proposition}
  \label[proposition]{join-fiber-sequences}
  Given two fiber sequences
  \[
    \begin{tikzcd}
      A_1\ar[r]&B_1\ar[r,"f_1"]&C
    \end{tikzcd}
    \qquad\qquad
    \begin{tikzcd}
      A_2\ar[r]&B_2\ar[r,"f_2"]&C
    \end{tikzcd}
  \]
  their join
  \[
    \begin{tikzcd}
      A_1\join A_2\ar[r]&B_1\join_C B_2\ar[r,"f_1\join f_2"]&C
    \end{tikzcd}
  \]
  is a fiber sequence.
\end{proposition}

Now consider two maps $\phi_i:A_i\to\B G$ for $i\in\set{1,2}$ and some group~$G$: these maps can be understood as actions of $G$ on $F_i\defd\fib{\phi_i}$ by action-fibration duality. Their join $\phi_1\join\phi_2$ corresponds to the action of~$G$ on $F_1\join F_2$ obtained by letting~$G$ act as $\phi_i$ on $F_i$ seen as sitting inside of $F_1\join F_2$ through the canonical map $F_i\to F_1\join F_2$ (given by the definition as a HIT), and extended ``by continuity'' on other points.

Given a map $f:A\to B$, we write $f^n$ for the iterated join of~$n$ copies of~$f$. Combining \cref{join-connected,fiber-join}, we have that $f^n$ is more and more connected as $n$ increases.
Taking the inductive limit as $n$ goes to infinity, it can be shown that these maps converge toward the inclusion $\im(f)\into B$~\cite[Theorem 4.2.13]{rijke2018classifying}.

\subsection{Higher-dimensional hypercubical manifolds in type theory}
We define the $n$-th higher hypercubical manifold as the type
\[
  \HM[n]
  \defd
  \HM\join_{\B Q}\ldots\join_{\B Q}\HM
\]
as the iterated join of $n$ copies of~$\HM$, and
\[
  \phi^n\defd\phi\join\ldots\join\phi:\HM[n]\to\B Q
\]
as the corresponding universal map. By iteratively joining the fundamental fiber sequence of \cref{the-fiber-sequence} with itself using \cref{join-fiber-sequences}, we obtain:

\begin{proposition}
  We have, for every $n\in\N$, a fiber sequence
  \[
    \begin{tikzcd}
      \S{4n-1}\ar[r]&\HM[n]\ar[r,"\phi^n"]&\B Q
    \end{tikzcd}
  \]
\end{proposition}

\noindent
As a direct consequence of the previous fiber sequence, and the fact that $\S k$ is $(k{-}1)$-connected~\cite[Corollary 8.2.2]{hottbook}, we have

\begin{proposition}
  The canonical map $\phi^n:\HM[n]\to\B Q$ is $(4n{-}2)$-connected.
\end{proposition}

\noindent
By taking the inductive limit as explained above, we obtain a delooping $K^\infty$ of~$Q$:

\begin{proposition}
  \label[proposition]{join-equiv}
  We have a fiber sequence
  \[
    \begin{tikzcd}
      1\ar[r]&K^\infty\ar[r,"\phi^\infty"]&\B Q
    \end{tikzcd}
  \]
  and thus $\phi^\infty$ is an equivalence.
\end{proposition}
\begin{proof}
  The fact that $\phi^\infty$ is an equivalence follows from the fact that it merely has a point and $\B Q$ is connected.
\end{proof}


\begin{remark}
  Note that the previous construction would have worked starting from any map $\phi:A\to\B Q$ (as long as~$A$ is inhabited), the pointing map $1\to\B Q$ for instance, in the sense that \cref{join-equiv} would still hold.
\end{remark}

\section{Conclusion}
\label{sec:conclusion}
We have defined the hypercubical manifold in homotopy type theory and shown the relevance of our construction by exhibiting the resulting type $\HM$ as a quotient of the sphere $\S3$ under the canonical action of the fundamental group~$Q$. We have also introduced higher-dimensional variants of these spaces.

As a direct application, it would be interesting to use $\HM[n]$ in order to compute the cohomology groups of $Q$ in low dimensions (up to $4n-2$). Indeed, from the definition of the delooping $\B G$ of any group~$G$, one can obtain a cellular type (the analogue, in type theory, of a CW-complex) by applying the Milnor construction~\citep{buchholtz2017real}. In our case, this provides a cellular description of~$\B Q$ from which one can synthetically compute, in homotopy type theory, the cellular cohomology of~$Q$, as the cohomology of $\B Q$ by applying the constructions of~\cite{buchholtz2020cellular}. With integral coefficients, all these groups are already known, and we should for instance be able to recover $H^2(Q)$. This would first require showing that the types~$\HM[n]$ are cellular, which is expected to hold, but is left for future work.

We believe that the methodology introduced here is very general and should be applicable in order to define types corresponding to various well-known spaces, which can be obtained from polyhedra by gluing faces, including those considered by~\cite{poincare1895analysis}. In particular, we plan to investigate a definition of the homology sphere~\citep{poincare1904cinquieme}, which is defined in a similar way, although the combinatorics is more involved: it can be described as a cellular complex with 5 0-cells, 10 1-cells, 6 2-cells and one 3-cell, its fundamental group has 120 elements and the fiber of the fundamental fibration is a model of~$\S3$ with 600 0-cells, 1200 1-cells, 720 2-cells and 120 3-cells. Another family of examples is given by suspensions of linear morphisms of the torus~\citep{arnold1967problemes}.

This paper also fits into the general program aiming at constructing the group $\SU(2)$ in homotopy type theory, along with its higher coherences. This amounts to constructing a delooping of the 3-sphere, which is known to be isomorphic to $\SU(2)$, and would thus have a group structure induced by its identification as a loop space. Constructing such a delooping is thus an important open problem in homotopy type theory. It is also equivalent to defining all quaternionic projective spaces via the Milnor construction, by generalizing the construction performed in the real case~\citep{buchholtz2017real}. A first important step in this direction is the construction, and classification, of the H-space structures on~$\S3$~\citep{buchholtz2016cayley,buchholtz2025central}. The group $\SU(2)$ canonically acts on itself (by left multiplication) resulting in a canonical action of $\SU(2)$ on~$\S3$, of which the action of~$Q$ on~$\S3$ is a special case (as are the actions of~$\Z_n$ on~$\S3$ constructed by~\cite{lens}). In particular, if we knew how to deloop $\S3$, we could hope for a coherent version of the action of~$\SU(2)$ on itself, from which one could deduce all the free actions of finite groups on~$\S3$ (because these are known to be subgroups of~$\SU(2)$~\citep{thurston2000three}).

\bibliographystyle{plainurl}
\bibliography{papers}
\end{document}